\documentclass[12pt]{article}

  \usepackage{amssymb}         
  \usepackage{amsmath}          
  \usepackage{amsfonts}           
  \usepackage{amsthm}
  \usepackage{amscd}
  \usepackage{enumerate}
  \usepackage[active]{srcltx}

  \usepackage[mathscr]{eucal}

  \usepackage{graphicx,epsfig}

\def\Teichmuller{{Teichm{\"u}ller}}
\def\residuiteratif{{r{\'e}sidue it{\'e}ratif}}

\def\mult{{\rm mult}}
\def\index{{\rm index}}
\def\resit{{\rm {r\acute{e}sit}}}

\def\Res{{\rm {Res}}}

  \newcommand*{\al}{\alpha}

  \newcommand{\Arg}{{\operatorname{Arg}}}

  \newcommand*{\be}{\beta}

  \newcommand*{\chdel}{\ensuremath{{\check{\delta}}}}
  \newcommand*{\chga}{\ensuremath{{\check{\gamma}}}}
  \newcommand*{\chphi}{\ensuremath{{\check{\phi}}}}
  \newcommand*{\chpsi}{\ensuremath{{\check{\psi}}}}
  \newcommand*{\cint}{\oint}
  \newcommand*{\C}{\mathbb{C}}
  \newcommand*{\Chat}{\ensuremath{{\widehat{\mathbb  C}}}}

  \renewcommand*{\d}{{\operatorname{d}}}
  \newcommand*{\diam}{{\operatorname{diam}}}

  \newcommand*{\dw}{{dw}}
  
  \newcommand*{\dz}{{dz}}
  \newcommand*{\D}{\mathbb{D}}
  \newcommand*{\Dbar}{{\overline{\mathbb{D}}}}

  \newcommand*{\Dstar}{{\mathbb{D}^*}}

  \newcommand*{\e}{\ensuremath{{\operatorname{e}}}}
  \newcommand*{\eps}{{\epsilon}}

  \newcommand*{\fla}{{f_\la}}

  \newcommand*{\ga}{\gamma}

  \newcommand*{\Ga}{\Gamma}

  \renewcommand*{\H}{\ensuremath{\mathbb  H}}
  \newcommand*{\Hminus}{\ensuremath{{{\mathbb  H}_-}}}
  \newcommand*{\Hplus}{\ensuremath{{{\mathbb  H}_r}}}
  \newcommand*{\Hpq}{{H_{p/q}}}

  \newcommand*{\Klap}{{K_\la'}}

  \newcommand*{\la}{\lambda}

  \newcommand*{\lastar}{{\lambda^*}}
  \newcommand*{\La}{\Lambda}

  \newcommand*{\Lpq}{{\mathrm{L}_{p/q}}}
  \newcommand*{\Lpqppqp}{{\mathrm{L}}_{p'/q'}^{p/q}}
  \newcommand*{\Log}{{\operatorname{Log}}}
  \newcommand*{\Lastar}{{\Lambda^*}}

  \newcommand*{\mapfromto}[3]{\hbox{\ensuremath{#1 : #2 \longrightarrow #3}}}

  \newcommand*{\mustar}{{\mu^*}}
  \newcommand*{\Mbrot}{{\mathrm{M}}}
  \newcommand*{\MHp}{{\Mbrot'_H}}
  \newcommand*{\ML}{{\Mbrot^L}}
  \newcommand*{\Mp}{{\Mbrot'}}
  \newcommand*{\Mpp}{{\Mbrot''}}
  
  \newcommand*{\Mod}{{\operatorname{Mod}}}

  \newcommand*{\Mobius}{{M{\"o}bius}}
  \newcommand*{\Mpq}{{\Mbrot_{p/q}}}
  \newcommand*{\Mppq}{{\Mbrot_{p'/q}}}
  \newcommand*{\Mstar}{{M^*}}

  \newcommand*{\MPQ}{{\Mbrot_{P/Q}}}
  \newcommand*{\Mz}{{\Mbrot_{0}}}
  \newcommand*{\N}{{\mathbb{N}}}
  
  \newcommand*{\om}{{\omega}}
  
  \newcommand*{\opq}{{\omega_{p/q}}}

  \newcommand*{\Om}{\Omega}
  \newcommand*{\OO}{{\mathcal{O}}}

  \newcommand*{\psipq}{{\psi_{p/q}}}
  \newcommand*{\psiPQ}{{\psi_{P/Q}}}
  
  \renewcommand*{\P}{{\mathrm{P}}}

  \newcommand*{\Pla}{{\mathrm{P}_\la}}
  \newcommand*{\Plaq}{{\mathrm{P}_\la^q}}
  \newcommand*{\Pmu}{{\mathrm{P}_\mu}}

  \newcommand*{\PPP}{{\mathcal{P}}}
  \newcommand*{\PSL}{{\operatorname{PSL}}}

  \newcommand*{\R}{\mathbb{R}}
  
  \newcommand*{\Rminus}{\ensuremath{{\mathbb  R}_-}}
  \newcommand*{\Rplus}{\ensuremath{{\mathbb  R}_+}}

  \newcommand*{\si}{\sigma}
  \newcommand*{\sm}{\ensuremath{{\smallsetminus}}}
  
  \newcommand*{\Sm}{\ensuremath{{\setminus}}}

  \newcommand*{\T}{{\mathbb{T}}}

  \newcommand*{\whal}{\ensuremath{{\widehat{\alpha}}}}

  \newcommand*{\whga}{\ensuremath{{\widehat{\ga}}}}
  \newcommand*{\whGa}{\ensuremath{{\widehat{\Ga}}}}

  \newcommand*{\whHpm}{\ensuremath{{\widehat{H}^\pm}}}
  \newcommand*{\whHpq}{{\widehat{H}_{p/q}}}

  \newcommand*{\whLa}{\ensuremath{{\widehat{\Lambda}}}}
  \newcommand*{\whLpq}{{\widehat{\mathrm{L}}_{p/q}}}
  \newcommand*{\whLpqn}{{\widehat{\mathrm{L}}_{\frac{n^2-1}{n^3}}^{p/q}}}
  \newcommand*{\whLPQn}{{\widehat{\mathrm{L}}_{\frac{n^2-1}{n^3}}^{P/Q}}}
  \newcommand*{\whLpqppqp}{{\widehat{\mathrm{L}}_{p'/q'}^{p/q}}}
  \newcommand*{\whMpq}{{\widehat{\Mbrot}_{p/q}}}
  \newcommand*{\whMPQ}{{\widehat{\Mbrot}_{P/Q}}}

  \newcommand*{\whxi}{\ensuremath{{\widehat{\xi}}}}
  
  \newcommand*{\wtga}{\ensuremath{{\widetilde{\ga}}}}
  \newcommand*{\wtGa}{\ensuremath{{\widetilde{\Ga}}}}
  \newcommand*{\wtphi}{\ensuremath{{\widetilde{\phi}}}}
  \newcommand*{\wtvarphi}{\ensuremath{{\widetilde{\varphi}}}}
  \newcommand*{\wtsi}{\ensuremath{{\widetilde{\si}}}}

  \newcommand*{\WpqMz}{{\mathrm{W}_{p/q}^\Mz}}
  
  \newcommand*{\Z}{\mathbb{Z}}

\newtheorem{definition}{Definition}
\newtheorem{theorem}{Theorem}
\newtheorem{corollary}[theorem]{Corollary}
\newtheorem{lemma}{Lemma}
\newtheorem{proposition}[theorem]{Proposition}

\newtheorem{remark}{Remark}
\newtheorem{conjecture}{Conjecture}

\newcommand{\ALIGN}{\begin{align*}}
\newcommand{\ENDALIGN}{\end{align*}}
\newcommand{\ENUM}{\begin{enumerate}}
\newcommand{\ENUMa}{\begin{enumerate}[a.]}
\newcommand{\ENUMA}{\begin{enumerate}[A.]}
\newcommand{\ENUMAF}{\begin{enumerate}[\bf A.]}
\newcommand{\ENUMi}{\begin{enumerate}[i)]}
\newcommand{\ENDENUM}{\end{enumerate}}
\newcommand{\ITMZ}{\begin{itemize}}
\newcommand{\ENDITMZ}{\end{itemize}}
\newcommand{\REFEQN}[1] { \begin{equation}\label{#1} }
\newcommand{\ENDEQN}{\end{equation}}
\newcommand{\THM}{\begin{theorem}}
\newcommand{\REFEXA}[1] { \begin{example}\label{#1} }
\newcommand{\ENDEXA}{\end{example}}
\newcommand{\MTX}{ \begin{matrix}}
\newcommand{\ENDMTX}{ \end{matrix}}
\newcommand{\REM}{ \begin{remark}}
\newcommand{\ENDREM}{\end{remark}}
\newcommand{\REFTHM}[1] { \begin{theorem}\label{#1} }
\newcommand{\RREFTHM}[2] { \begin{theorem}[#1]\label{#2} }
\newcommand{\ENDTHM}{\end{theorem}}
\newcommand{\REFNTH}[1] { \begin{newthm}\label{#1} }
\newcommand{\ENDNTH}{\end{newthm}}
\newcommand{\REFPROP}[1]{\begin{proposition}\label{#1} }
\newcommand{\RREFPROP}[2]{\begin{proposition}[#1]\label{#2} }
\newcommand{\PROP}{\begin{proposition}}
\newcommand{\ENDPROP}{\end{proposition} }
\newcommand{\REFDEF}[1]{\begin{definition}\label{#1} }
\newcommand{\DEF}{\begin{definition}}
\newcommand{\ENDDEF}{\end{definition} }
\newcommand{\REFLEM}[1]{\begin{lemma}\label{#1} }
\newcommand{\RREFLEM}[2]{\begin{lemma}[#1]\label{#2} }
\newcommand{\LEM}{\begin{lemma}}
\newcommand{\ENDLEM}{\end{lemma} }
\newcommand{\REFCOR}[1]{\begin{corollary}\label{#1} }
\newcommand{\COR}{\begin{corollary}}
\newcommand{\ENDCOR}{\end{corollary} }
\newcommand{\CONJ}{\begin{conjecture}}
\newcommand{\REFCONJ}[1]{\begin{conjecture}\label{#1}}
\newcommand{\RREFCONJ}[2]{\begin{conjecture}{#1}\label{#2}}
\newcommand{\ENDCONJ}{\end{conjecture} }
\newcommand{\REFDEFTHM}[1] { \begin{defthm}\label{#1} }
\newcommand{\ENDDEFTHM}{\end{defthm}}

\newcommand{\corref}[1]{Corollary~\ref{#1}}
\newcommand{\conjref}[1]{Conjecture~\ref{#1}}

\newcommand{\thmref}[1]{Theorem~\ref{#1}}
\newcommand{\propref}[1]{Proposition~\ref{#1}}

\newcommand{\PROOF}{\begin{proof}}
\newcommand{\ENDPROOF}{\end{proof}}

\newcommand{\Bottcher}{B{\"o}ttcher}
\usepackage{color} 
\title{On quasi-conformal (in-) compatibility
of satellite copies of the Mandelbrot set: I} \author{Luna Lomonaco\qquad Carsten Lunde Petersen}

\begin{document}

\maketitle

\begin{abstract}
In the paper \cite{DouadyandHubbard} Douady and Hubbard introduced the notion of polynomial-like maps. 
They used it to identify homeomophic copies $\Mp$ of the Mandelbrot set 
inside the Mandelbrot set $\Mbrot$. 
They conjectured that in case of primitive copies the homeomorphism 
between $\Mp$ and $\Mbrot$ is q-c, 
and similarly in the satellite case, it is q-c~off any small neighborhood of the root.
These conjectures are now Theorems due to Lyubich, \cite{Lyubich}.
The satellite copies $\Mpq$ are clearly not q-c homeomorphic to $\Mbrot$. 
But are they mutually q-c homeomorphic? 
Or even q-c homeomorphic to half of the logistic Mandelbrot set? 
In this paper we prove that, in general, 
the induced Douady-Hubbard homeomorphism is not the restriction of a q-c homeomorphism: 
For any two copies $\Mp$ and $\Mpp$ the induced Douady-Hubbard homeomorphism 
is not q-c, if the root multipliers $\la'$ and $\la''$, 
which are primitive $q$ and $q'$ roots of unity, have $q\not= q'$.
\end{abstract}

\setcounter{tocdepth}{1}
\tableofcontents

\section{Introduction}\label{intro}
We consider the moduli space of quadratic polynomials as dynamical systems under iteration. 
This moduli space is conveniently parametrized by the family $Q_c(z) = z^2+c$ 
of monic quadratic polynomials with critical point $0$ and critical value $c\in\C$.
For an introduction to the theory of holomorphic dynamical systems given by iteration 
see e.g. \cite{CarlesonandGamelin}, \cite{Milnor}.
The connectedness locus of the quadratic family $Q_c$, $c\in\C$ is the 
famous Mandelbrot set:
$$
\Mbrot=\{c\in\C| J_c \textrm{ is connected}\}.
$$
Conjecturally the interior of $\Mbrot$ consists of precisely the hyperbolic components 
generically denoted $H$, 
each of which is conformally equivalent to $\D$ through the map which assigns to 
$c\in H$ the multiplier $\rho_H(c)$ of the unique attracting cycle for $Q_c$. 
The period $m\geq 1$ of this attracting cycle is a characteristic of $H$. 
The multiplier map is known to extend to a homeomorphism between the closures. 
Moreover this map is even a holomorphic diffeomorphism between the closures except possibly 
at the root point with image $1$, which may be either a cusp point, the so-called 
primitive case, or a holomorphically smooth point, the so-called satellite case. 
For a hyperbolic component $H$ of period $m$, each boundary point $\rho_H^{-1}(\opq)$ 
where $\opq= \e^{i2\pi p/q}$ is the root point 
of a satellite hyperbolic component $H'$ of $H$ with characteristic period $qm$.
The main hyperbolic component $H_0$ of $\Mbrot$ is the unique component with 
characteristic period $1$. It is the primitive component bounded by the central cardiod. 

Douady and Hubbard \cite{DouadyandHubbard} identified via the theory of renormalization 
for each hyperbolic component $H$ a subset denoted $\MHp$ and constructed 
on $\MHp$ an injection {\mapfromto {\chi_H} {\Mbrot_H}{\Mbrot}} 
with $\chi_H(H) = H_0$. 
They showed using the theory of polynomial-like mappings, 
that $\chi$ is a homeomorphism if $H$ is primitive and 
except possibly above the root of $H_0$, if $H$ is a satellite component. 
They moreover conjectured that $\chi_H$ is the restriction of a quasi-conformal map 
except possibly on neighborhoods of the root, when $H$ is a satellite. 
Haissinsky \cite{Haissinsky} complemented the above by showing that $\chi_H$ 
is also a homeomorphism at the root in the satellite case. 
And Lyubich \cite{Lyubich} proved the long standing conjecture 
on quasi-conformality of $\chi_H$. Anyway, Lyubich proof fundamentaly breaks down for 
neighborhoods of roots of satellite hyperbolic components.

The question which remains is whether the satellite copies $\Mpq$ with 
roots $\la_{H_0}^{-1}(\opq)$ are mutually quasi-conformally homeomorphic. In particular, do they have a common model with which they are 
quasi-conformally homeomorphic to? 
A natural such candidate is given by half the logistic Mandelbrot set.
Let $P_\la(z) = \la z + z^2$ denote the logistic family of quadratic polynomials, 
which is related to the family $Q_c$ by the degree $2$ branched covering 
$c(\la) = \la/2 - \la^2/4$. 
The connectedness locus is the logistic Mandelbrot set $\ML$, 
which naturally decomposes into a left and a right half 
with unique common point $1$ and each of which is mapped 
homeomorphically onto $\Mbrot$ by the holomorphic map $\la\mapsto c(\la)$. 

Using her newly developed notion of parabolic-like mappings, 
the first author have shown that 
the root dynamics of any two satellite copies, are q-c conjugate, \cite{Lomonaco}. 
Nevertheless we prove in this paper that:



\REFTHM{main}
For $p/q$ and $P/Q$ irreducible rationals with $q \not= Q$, 
the induced Douady-Hubbard homeomorphism 
$$
{\mapfromto {\xi = \chi_{P/Q}^{-1}\circ\chi_{p/q}} \Mpq \MPQ},$$
is not quasi-con\-formal, \textit{i.e.} it does not admit a quasi-conformal 
extension to any neighborhood of the root.
Moreover for any satellite copy the induced Douady-Hubbard homeomorphism 
to half the logistic Mandelbrot set is not a q-c homeomorphism. 
\ENDTHM

We have reasons which advocates the following 
\CONJ\label{reducedconj}
For $p/q$ and $p'/q$ irreducible rationals with the same denominator, 
the induced Douady-Hubbard homeomorphism 
$$
{\mapfromto {\xi = \chi_{p'/q}^{-1}\circ\chi_{p/q}} \Mpq \Mppq}
$$ 
does admit a quasi-conformal extension between neighborhoods 
of the satellite copies.
\ENDCONJ

We return to this conjecture in a second paper on the topic. 
The insight behind the proof, can be explained in brief as follows. 
Consider (as in the proof of Yoccoz inequality) the lift of the satellite copy $\Mpq$ 
to the $\log \la$ plane, i.e consider 
the connected component containing $i2\pi p/q$ 
of the pre-image of $\Mpq$ under the map $\La\mapsto \e^\La/2 + \e^{2\La}/4$. 
Rescale and relocate this pre-image by the affine map $\La\mapsto q\La-pi2\pi$, 
which is natural from the point of view of renormalization. 
Denote this lifted rescaled and translated copy with root $0$ by $\whMpq$. 
For $p/q$ and $P/Q$ irreducible rationals let {\mapfromto \whxi \whMpq \whMPQ}
denote the induced Douady-Hubbard homemorphism. 
Then a necessary condition for $\xi$ to be the restriction of a q.-c.~homeomorphism, 
is that $\whxi$ is uniformly close to the identity with respect to the hyperbolic metric 
on the right half-plane $\H_r:=\{x+iy| x>0\}$. 

As a word of further motivation we recall that Branner and Fagella \cite{BrannerandFagella} proved that, 
for $p/q$ and $p'/q$ irreducible rationals with the same denominator, the corresponding limbs $\Lpq$ and $L_{p'/q}$ of $H_0$ are homeomorphic 
and they conjectured their homeomorphisms are restrictions of quasi-conformal maps. 
\conjref{reducedconj} is a weaker, but yet essential version of the Branner-Fagella conjecture. \\

\textbf{Aknowledgment.} The authors would like to thank M. Lyubich 
for bringing the question this article answers under the attention of
the first author, when she was visiting Stony Brook University (fall 2013). 
The first author would also like to thank the Institute of Mathematical Sciences of Stony Brook University 
for their hospitality, and the Fapesp for support via the process 2013/20480-7.
The second author would like to thank 
the Danish Council for Independent Research | Natural Sciences for support via the 
grant 10-083122 and the Institute of Mathematical Sciences of Stony Brook  University 
for their hospitality. 

\section{Outline of the Proof of the main Theorem}
For the ease of exposition it is convenient to formulate and prove Theorem~\ref{main} 
in terms of the logistic or fixed point normalized quadratic family:
$$
\Pla(z)=\lambda z+z^2, \la\in\C
$$ 
where the parameter $\la$ is the multiplier of the fixed point $0$. 
This family forms a branched double covering of the standard representation $Q_c$, 
the relation being $c(\la) = \frac{\la}{2}-\frac{\la^2}{4}$, with critical point $1$ 
and critical value $1/4$. 
We write \textit{$\Mbrot_0$} for the homeomorphic pre-image 
$\la^{-1}(\Mbrot) \subset \{\la| \Re(\la)\leq 1\}$, and 
denote the other copy $\Mbrot_{0/1} := \la^{-1}(\Mbrot) \subset \{\la| \Re(\la)\geq 1\}$. 
Evidently $H_0=c(\D)$.
For $\la\in\Mbrot_{0/1}$ the fixed point $0$ 
is the $\beta$ fixed point of $P_\la$, that is the fixed point, 
which is the landing point of the unique fixed ray. 
For $\la\in\Mbrot_0$ the fixed point $0$ is the \emph{$\alpha$-fixed point} for $\Pla$, 
i.e.~the other fixed point. 

We shall abuse the notation and for each irreducible rational $p/q$ 
write $\Mpq \subseteq \Lpq$ also for the satellite copy 
and limb of $\Mz$ with root $\opq$. 
The limb $\Lpq$ is the unique connected component 
of $\Mz\Sm\{\opq\}$ not containing $0$. 
 Note that for $p/q=0/1$ we have $\Mpq = \Lpq$, and otherwise the inclusion is strict.

For $p/q$ and $P/Q$ two distinct irreducible rationals we likewise denote by $\xi$ 
the induced map {\mapfromto {\xi = \chi_{P/Q}^{-1}\circ\chi_{p/q}} \Mpq \MPQ} 
between the new sets. 
Since $\la\mapsto c(\la)$ is itself a quadratic polynomial, 
which maps the half-plane $\{x+iy| x<1\}$ univalently onto $\C\Sm[\frac{1}{4},\infty[$,
we may rephrase Theorem \ref{main} as:
\REFTHM{enlmain}
For $p/q$ and $P/Q$ irreducible rationals with $q \not= Q$ 
the Doaudy-Hubbard homeomorphism 
{\mapfromto {\xi = \chi_{P/Q}^{-1}\circ\chi_{p/q}} \Mpq \MPQ},
is not quasi-con\-formal,  \textit{i.e.} it does not admit a quasi-conformal 
extension between any neighborhoods of the roots.
\ENDTHM

\vspace{0.5cm}

For $p/q$ an irreducible rational and $\lambda \in \Mpq$, 
let {\mapfromto {f=f_\la} {U'} U} be a polynomial-like restriction of $\Plaq$ 
with the critical value of $\Pla$ in $U'$. 
Then the $\beta$-fixed point of the polynomial-like restriction is $0$, 
(the $\al$-fixed point of $\Pla$, when $p/q\not=0/1$) and its new multiplier is $\la^q$. 
Define $\La(\la)=\Log(\la^q)$ as the principal logarithm of the new multiplier. 
It follows from Yoccoz inequality that the map $\La$ is univalent on a neighborhood 
of the limb $\Lpq$. 
We define $\whMpq := \La(\Mpq)\subset\whLpq := \La(\Lpq)\subset\Hplus\cup\{0\}$ 
and may henceforth use $\La$ at will as the parameter.

For $P/Q$ another irreducible rational 
we shall, in order to avoid heavy indexing, 
use $\mu$ rather than $\la$ as the parameter 
and for $\mu\in\MPQ$ write {\mapfromto {g=g_\mu} {V'} V} 
for the quadratic like restriction 
and $M(\mu) = \Log(\mu^Q)$ as above and $\whMPQ = M(\MPQ)$.
Moreover we shall write {\mapfromto \whxi \whMpq \whMPQ} for 
the induced Douady-Hubbard homeomorphism. 
Behind the scene of the proof of \thmref{main} 
hides the following Proposition, 
which is an instance of a Theorem by Sullivan and Shub, 
and for which the proof essentially is \Teichmuller's extremal theorem 
for complex tori (see \propref{Teichmuller}).
\REFPROP{lowerdilationbound}
Let $\La\in\whMpq\Sm\{0\}$ and $M=\whxi(\La) \in \whMPQ$.
Then for any quasi-conformal (in particular any hybrid) conjugacy 
{\mapfromto \phi U V} between polynomial-like 
restrictions $f$ and $g$ as above we have:
\REFEQN{sullivanbound}
\limsup_{z \rightarrow 0} \Log K_{\phi}(z) \geq d_{\Hplus}(\La, M),
\ENDEQN
 where $K_\phi$ denotes the real dilation of $\phi$.
\ENDPROP
\noindent 

We next obtain a similar result in parameter space at Misiurewicz parameters 
$\lastar$ for which the critical value of $f$ is pre-fixed to $0$. 
For the proof of this we first invoke a rather standard holomorphic motion argument 
to pass from a neighborhood of $\lastar$ in parameterspace to 
a neighborhood of $0$ in the linearizing coordinate for $f_\lastar$ at $\be':=0$. 
Next we combine this with a new theorem, which is a relative of the {\Teichmuller} 
extremal theorem for complex tori to obtain: 
(for details see page \pageref{proofofdyn}):
\REFTHM{dynamicsimplyparameter}
Let $\Lastar\in\whMpq$ be any (Misiurewicz) parameter such that 
the critical value $-\lastar^2/4$ for $f_\lastar$ is prefixed to $0$, and 
let $\Mstar = \whxi(\Lastar)$. Then :
$$
\limsup_{\La\rightarrow\Lastar} \Log K_{\whxi}(\La) \geq d_{\Hplus}(\Lastar,\Mstar).
$$
\ENDTHM

For $p/q$ any irreducible rational and $\la\in\Mpq$ 
let $\alpha'(\la)$ denote the $\al$-fixed point $f_\la$. 
The multiplier map {\mapfromto \rho \Hpq \D}, $\rho(\la) = f_\la'(\alpha'(\la))$ 
univalently extends to a neighborhood of the root $\lambda_0= e^{2 \pi i p/q}$.
Thus the inverse $\la(\rho)$ provides a local, univalent parametrization 
$\rho\mapsto \La(\rho) = \Log(\la^q(\rho))$ of a neighborhood of $0$ 
in the $\La$ parameter plane 
by a neighborhood $\Om$ of $1$ in the $\rho$ plane.
We have the following estimate for $\La(\rho)$:

\REFPROP{goodassymtoticestimates}
For any irreducible rational $p/q$ 
there exists a constant $\Res_{p/q}$ with 
positive real part such that:
$$
\La(\rho)= -\frac{\Log(\rho)}{q}
-\Res_{p/q}\frac{(\Log\rho)^2}{q^2} + 
O\left(\left(\frac{\Log\rho}{q}\right)^3\right).
$$
\ENDPROP
\noindent
For $p/q$ and $P/Q$ two distinct irreducible rationals 
we have local parametrizations $\La(\rho)$ and $M(\rho)$ of a neighborhood of $0$ 
in the $\La$-parameter plane and the $M$-parameter plane by a neighborhood 
of $1$ in the $\rho$-plane.
With $\whxi$ as above we have $M(\rho) = \whxi(\La(\rho))$ for any $\rho\in\Dbar$. 
By elementary calculations the above asymptotic formulas yield:

\REFCOR{multiplierhyperbolicdistancediverge}
 For $q\not= Q$ and $\rho=e^{it}\in\partial\D$: 
 $$
 d_{\Hplus }(\Lambda(\rho),M(\rho))\longrightarrow \infty 
 \quad\textrm{as}\quad \rho \rightarrow 1\qquad (\textrm{i.e.,~as}\quad \R\ni t\to 0).
 $$ 
\ENDCOR

Recall that for $p/q$ an irreducible rational the set $\Lpq\supset\Mpq$ is 
the $p/q$-limb of $\Mbrot_0$. 
For $p'/q'$ a second irreducible rational we denote by $\Lpqppqp$ 
the $p'/q'$-sublimb of $\Lpq$, i.e.~the connected component 
of $\Mz\Sm\{\la(\rho_{p'/q'})\}$ not containing $0$, 
where $\rho_{p'/q'}=  e^{2 \pi i p'/q'}$. 
As above we write $\whLpqppqp := \La(\Lpqppqp)$. 
The following Proposition, which is an elementary consequence of 
the Yoccoz parameter inequality for the quadratic family, 
shows that, at least for particular sequences of rationals converging to $0$, 
the limb $\whLpqppqp$ is contained in a small hyperbolic neighborhood of its root.

\REFPROP{hyperbolicallysmalllimbs} 
For any irreducible rational $p/q$ the $\frac{n^2-1}{n^3}$-sublimbs of $\Lpq$ 
satisfy
$$
\diam_\Hplus(\whLpqn),~
\longrightarrow~ 0\qquad\textrm{as}\quad n\rightarrow \infty.
$$
\ENDPROP

For $p/q$ and $P/Q$ irreducible rationals and {\mapfromto \whxi \whMpq \whMPQ} 
we obtain, by combining 
\propref{hyperbolicallysmalllimbs} with \corref{multiplierhyperbolicdistancediverge} :
\REFCOR{Dynamicaldivergenceinsmalllimbs}\label{cor}
Let ${\{\La_n\}}_n$ with $\La_n \in\whLpqn\cap\whMpq$ be arbitrary 
and let $M_n = \whxi(\La_n)\in\whLPQn\cap\whMPQ$. 
If $q\not= Q$ then 
$$
d_{\Hplus}(\La_n,M_n)\longrightarrow\infty\qquad\textrm{as}\qquad n \rightarrow \infty
$$ 
\ENDCOR

\thmref{enlmain} and thus \thmref{main} now follows from 
\propref{dynamicsimplyparameter} combined with \corref{cor} applied to any 
sequence $\La_n^*\in \whLpqn\cap\whMpq$ of Misiurewicz parameters 
such that the critical point is prefixed to the $\beta'$ fixed point $0$.

\section{Proofs}
\subsection{Notation and {\`a} priori classical results.}
For $P$ a monic polynomial of degree $d\geq 2$ let $J_P$ and $K_P$ 
denote the Julia set and the filled Julia set of $P$ respectively. 
Denote by $\varphi_P$ with $\varphi_P\circ P = (\varphi_P)^d$
the {\Bottcher} coordinate for $P$ 
tangent to the identity at $\infty$ and let $g_P$ denote the Greens function for $K_P$. 
For $\theta\in\T:=\R/\Z$ denote by $R_\theta = R_\theta^P$ 
the external ray of argument $\theta$. 
That is $R_\theta$ is the field-line for $g_P$ 
arriving to $\infty$ at angle $\theta$,
 the latter being equivalent to
$\varphi_P(R_\theta) \subset \{re^{i2\pi\theta}|r>1\}$. 
For $K_P$ connected the {\Bottcher} coordinate extends to an isomorphism 
{\mapfromto {\varphi_P} {\Chat\Sm K_P} {\Chat\Sm\Dbar}} 
and $R_\theta = \varphi_P^{-1}(\{re^{i2\pi\theta}|r>1\})$. 
The ray is said to land at $z_0$ if its accumulation set in $K_P$ is the singleton $\{z_0\}$.
The mapping $z\mapsto z^d$ induces a dynamics on rays given by $\theta\mapsto d\theta$.

The dynamical plane landing theorem of Sullivan and Douady-Hubbard 
asserts that, when $K_P$ is connected, 
any periodic ray lands at a repelling or parabolic periodic point 
with period dividing that of the ray, 
and any pre periodic ray lands at a point pre periodic to a repelling or parabolic point. 

A repelling or parabolic periodic point $z_0$ of exact period $k\geq 1$ 
for a polynomial $P$ is said to have combinatorial rotation number $p/q$, 
$(p,q)=1$, if it is the common landing point for a $q$ cycle of rays under $P^k$, and if the rays under $P^k$
are permuted 
 as the rational rotation $p/q$ 
in the counter clockwise order around $z_0$. 
The Douady landing theorem asserts that any (pre)-periodic point is the landing point 
of at least one and at most finitely many (pre) periodic rays. 
And this defines a combinatorial rotation number for $z_0$. 
For angle doubling ($z\mapsto z^2$) there is for each irreducible rational 
$p/q$ a unique cycle of arguments (of rays) with rational rotation number $p/q$.

For the quadratic polynomials $Q_c$ we write $\phi_c$ for $\phi_{Q_c}$ 
and $R_\theta^c$ for the external ray of argument $\theta$. 
The uniformization of the complement of the Mandelbrot set 
{\mapfromto \Phi {\Chat\Sm\Mbrot} {\Chat\Sm\Dbar}} is given by 
$\Phi(c) = \varphi_c(c)$. The external ray for $\Mbrot$ of angle $\theta$ 
is $R_\theta^\Mbrot := \Phi^{-1}(R_\theta^0)$.

The Doaudy-Hubbard landing theorem for $\Mbrot$ asserts that any periodic ray 
(ray of periodic argument) for $\Mbrot$ lands at the root of a hyperbolic component 
with the same period, and except for the special case of $\theta=0=1$, 
the root $c_H$ of a hyperbolic component $H$ of period $m$ 
is the landing point of a pair of parameter rays $R_\theta^\Mbrot$ and $R_{\theta'}^\Mbrot$. 
These rays separate the plane, and so they define a wake $W_{c_H}^\Mbrot$ 
containing $H$ but not $0$. Moreover, for $c_H$ 
the corresponding dynamical rays define a wake $W_{c_H}^{c_H}$ 
in which the critical value $c_H$ is the unique point 
among the first $m$ iterates under $Q_{c_H}$ of $0$. 
In addition the common parabolic landing point $z_0^H$ of the dynamical rays 
$R_\theta^{c_H}$, $R_{\theta'}^{c_H}$ has period $k$ dividing $m$. 
If $k=m$ the hyperbolic component is primitive. 
Else $H$ is a satellite component of a hyperbolic component $H'$ 
of period $k$ and attached at an internal angle $r/s$ with $sk=m$.
In this satellite case the set $\mathrm{L}_{c_H}^{H'}:=\Mbrot\cap W_{c_H}^\Mbrot$ 
is called the limb of $H'$ with root $c_H$.
 
For $c\in W_{c_H}^\Mbrot$ the rays $R_\theta^c$, $R_{\theta'}^c$ 
co-land on a holomorphically moving repelling and $k$-periodic point $z(c)$, 
which converges to $z_0^H$ as $c\to c_H$ in $W_{c_H}^\Mbrot$. 
The dynamical wake $W_{c_H}^c$ defined by the co-landing pair of rays 
$R_\theta^c$, $R_{\theta'}^c$ contain an $m$-periodic point $w(c)$, 
which likewise depends holomorphically on $c\in W_{c_H}^\Mbrot$, 
which likewise converges to $z_0^H$ as $c\to c_H$ in $W_{c_H}^\Mbrot$ 
and which is non repelling precisely for $c\in\overline{H}$. 
The multiplier map {\mapfromto {\rho_H} H \D}, $\rho_H(c) = (Q_c^m)'(w(c))$ 
thus extends to a holomorphic function on $W_{c_H}^\Mbrot$ and even 
to a neighborhood of the root $c_H$.

The Doaudy-Hubbard landing theorem for $\Mbrot$ also asserts that 
any pre-periodic parameter ray $R_\theta^\Mbrot$ of pre-period $l$ and period $k$ 
lands at a Misiurewich parameter $c$, such that $c$ is  
the landing point of the dynamical ray $R_\theta^c$ and it has 
pre period $l$ to a repelling periodic cycle of period dividing $k$. 
Moreover for any dynamical ray $R_{\theta'}^c$ landing at $c$ the corresponding 
parameter ray $R_{\theta'}^\Mbrot$ also lands at $c$. 

For $P$ a polynomial of degree $d\geq 2$ and $z_0$ a repelling periodic point with 
exact period $k$ and combinatorial rotation number $p/q$, 
the Yoccoz multiplier theorem or the PLY-inequality asserts that:
There exists a unique \emph{preferred logarithm} $\La$ 
of the multiplier $\la := (P^k)'(z_0)$ 
such that any lift of a ray landing on $z_0$ to the log linearizer $\Psi$ 
for $P^k$ at $z_0$ (i.e.~$P^k(\Psi(z)) = \Psi(z+\La)$ is invariant 
under the translation $z\mapsto z + q\La - i2\pi$. Moreover, this preferred logarithm is 
contained in the closed disk of radius $r$ and center $r+i2\pi p/q$, where 
\REFEQN{yoccozdynamical}
r = \frac{k\log(d)}{nq}
\ENDEQN
and $n$ is the number of cycles of external rays under $P^k$ co-landing on $z_0$ 
(see \cite{Petersen} and \cite{Hubbard} for more details). 

The Yoccoz parameter inequality, which immediately follows from the above, asserts 
that for any hyperbolic component $H$ of period $k$ 
there exists a constant $C=C_{H}>0$ such that for any $p/q$, $(p,q)=1$, 
the limb $\mathrm{L}_{p/q}^{H}$ of $H$ with internal angle $p/q$ satisfies
\REFEQN{yoccozparameterineq}
\diam_E(\mathrm{L}_{p/q}^{H}) \leq \frac{C_H}{q}.
\ENDEQN
For a proof, see e.g. \cite[Remark 4.3]{Hubbard} (in this proof 
there is an impreciseness in the bound in the primitive case, where the exponent should 
have been $2$ rather than $1/2$).

The notion of rays, limbs, etc.,~as well as Yoccoz parameter inequality (with new constant) 
naturally carries over to the logistic Mandelbrot set, $\Mz$, to $\log(\Mz)$ 
and to the rescaled and shifted log-limbs and log-$\Mpq$'s 
$\whMpq\subset \whLpq$.

It immediately follows from \eqref{yoccozdynamical} 
that the rescaled and shifted log-limbs and log-$\Mpq$'s 
$\whMpq\subset \whLpq$ are all contained in the closed disk of radius and 
center $r_0 = \log 2$. 
It follows from this that for $\la\in\Lpq$ the preferred logarithm of the 
multiplier $\la^q$ for $P_\la(z) = \la z + z^2$, is the principal logarithm 
$\La=\Log(\la^q)$.

\subsection{Proof of \propref{lowerdilationbound}}
Suppose {\mapfromto f U \C} is a holomorphic function with a repelling fixed point 
$z_0 \in U$ of multiplier $\la$. Suppose further that {\mapfromto \ga {[0,1[} \C} 
is an arc with $\ga(0) = z_0$ and $f(\ga([0,a[)) = [1,0[$ for some $a>0$. 
Then \emph{$\ga$ defines a preferred logarithm $\La$ of $\la$} as follows: 
The quotient space $T = T_f := \Dstar_r(z_0)/f$ of the punctured disk of center $z_0$ 
and sufficiently small radius $r>0$ by the dynamics of $f$ 
is a complex torus. 
The quotient $\ga/f$ is a simple closed curve $\chga\subset T$, 
which together with the projection $\chdel$ of a small clockwise oriented circle 
$C_{r'}(z_0)$, $0<r'<r$, generates the fundamental group of $T$. 
The triple $(T,\chdel, \chga)$ defines a point $\La\in \H_r$, 
in the {\Teichmuller} space of complex tori modeled on $\H_r$ by taking 
$-i2\pi$ and $\La'\in\H_r$ as generators. 
The preferred logarithm of $\la$ defined by $\ga$ is the {\Teichmuller} point
$\La$ represented by  $(T,\chdel, \chga)$. 

In more detail, possibly reparametrizing $\ga$ we can suppose 
$\ga(\e t) = f\circ\ga(t)$. 
Let {\mapfromto \psi {\D} U} with $\psi(0)=z_0$ 
be a univalent local linearizer for $f$ around $z_0$, i.e. $\psi(\la z) = f\circ\psi(z)$. 
Let {\mapfromto \whga {[0,\eps]} {\D_r}} be the induced parametrization 
of the connected component of $\psi^{-1}(\ga)\cap\D_r$ containing $0$. 
Using the dynamics we can extend $\whga$ to 
an arc {\mapfromto \whga {[0,\infty]} {\Chat}} with $\whga(\e t) = \la\whga(t)$ and 
$\whga(\infty) = \infty$. 
For any logarithm $\La'$ of $\la$ the exponential map $\e^z$ semi-conjugates 
translation by $\La'$ to multiplication by $\la$. 
Let {\mapfromto \wtga \R \C} be a lift of $\whga(\e^t)$ to $\e^z$. 
Then all such lifts are disjoint and given by $\wtga_n := \wtga + n i2\pi$, 
$n \in \Z$ and there exists $m = m(\La')\in\Z$ such that 
$\La' + \wtga_n(t) = \wtga_{n+m}(t+1)$ for all $n$. 
Thus replacing $\La'$ by $\La = \La' - mi2\pi$ we obtain 
$\La + \wtga_n(t) = \wtga_n(t+1)$ for all $n$. 
\emph{The preferred logarithm of $\la$ defined by $\ga$ is $\La$.} 
Let $\Ga := \{ m\La - n2\pi i | m, n \in \Z \}$ and define the complex torus
$T_\La := \C/\Ga$, then $\La\R/\Ga$ and $\wtga/\Ga$ are isotopic in $T_\La$. 
Moreover the mapping {\mapfromto {\psi(\e^z)} {\Hminus} {\psi(\D)\Sm\{z_0\}}} 
induces an isomorphism {\mapfromto \chpsi {T_\La} {T_f}}, 
which maps $\wtga/\Ga$ onto $\chga$ and $[0, -i2\pi]/\Ga$ into an arc 
isotopic to $\chdel$. That is $(T,\chdel, \chga)$ represents the {\Teichmuller} point
$\La$.

Next suppose {\mapfromto {f_j} {U_j} \C}, $j= 1,2$ are holomorphic functions 
with respective repelling fixed points $z_j \in U_j$ of multiplier $\la_j$. 
Suppose moreover that $\La_J$ are preferred logarithms of $\la_j$ 
defined by backwards invariant arcs {\mapfromto {\ga_j} {[0,1[} \C} 
with $\ga_j(0) = z_j$. 
Then any quasi conformal local conjugacy {\mapfromto \phi {\om_1} {\om_2}} 
of $f_1$ to $f_2$ between neighborhoods $\om_1$ of $z_1$ and $\om_2$ of $z_2$
induces a quasi-conformal homeomorphism 
{\mapfromto {\chphi} {T_1} {T_2}} sending $\chdel_1\subset T_1 := T_{f_1}$ to an 
arc isotopic to $\chdel_2\subset T_2 := T_{f_2}$.

\REFPROP{Teichmuller}
With the notation above, if $\chphi(\chga_1)$ is isotopic to $\chga_2$ in $T_2$, or 
equivalently if one (hence any) lift $\wtphi$ of $\phi\circ\psi_1(\e^z)$ to $\psi_2(\e^z)$ 
conjugates translation by $\La_1$ to translation by $\La_2$,
then the real dilatation $K_\phi$ satisfies
$$
\limsup_{|z|\to 0} K_\phi(z) \geq \d_{\Hplus}(\La_1, \La_2),
$$
where $\d_\Hplus(\cdot,\cdot)$ denote the hyperbolic distance in $\Hplus$.
\ENDPROP
\PROOF
This is an immediate consequence of the {\Teichmuller} extremal Theorem for 
complex tori, which states that real dilatation $K_{\chphi}$ of $\chphi$ satisfies
$$
\sup_{z\in T_1} K_{\chphi}(z) \geq \d_{\Hplus}(\La_1, \La_2).
$$
\ENDPROOF

\PROOF{of \propref{lowerdilationbound} :}
Any ray in the $p/q$-cycle of rays for $P_\la$ landing at $0$ is 
fixed under $f_\la$ and defines $\La$ as the preferred logarithm of $\la^q$. 
Similarly $M$ is the preferred logarithm of $\mu^Q$ defined by 
the rays which colland on $0$ and which are fixed by $g_\mu$. 
Moreover any conjugacy between $f_\la$ and $g_\mu$ preserves 
the filled Julia sets and so the isotopy hypothesis of the above \propref{Teichmuller}
is satisfied. Hence \propref{lowerdilationbound} follows.
\ENDPROOF

\subsection{Proof of \thmref{dynamicsimplyparameter}}\label{Review}
Let $p/q$ be any irreducible rational and let $\la\in\Mpq$. 
Recall that $\al'_\la$ denotes the $\alpha$-fixed point of $f_\la$ 
and that $\al'(\la)$ can be followed holomorphically on the wake 
$W_{p/q}^\Mz$ with root $\la_0=\opq$ as a periodic point 
of exact period $q$ for $P_\la$. 
And that its multiplier $\rho(\la) := P_\la^q(\al'(\la))$ 
is a holomorphic function {\mapfromto \rho \WpqMz \C}, 
in fact with a holomorphic extension to a neighborhood of the root $\opq$,
with $\rho^{-1}(\Dbar) = \overline{\Hpq}$.

For any irreducible rational $p'/q'$ denote by $\Lpqppqp$ 
the $p'/q'$-sub-limb of $\Lpq=\Mz\cap\WpqMz$ 
with root $\la=\rho^{-1}(\om_{p'/q'})$.
For any $\la\in \Lpqppqp$ the $f_\la$ -fixed point $\al'(\la)$ 
is the common landing point of the $q'$-cycle of periodic rays for $f_\la$ 
($qq'$-cycle of rays for $\Pla$) giving $\al'(\la)$ rotation number $p'/q'$. 
In particular it is a separating point for the filled Julia sets $\Klap$ for $f_\la$ 
and $K_\la$ for $P_\la$. It separates $-(\la/2)^2$, 
the common critical value of $\Pla$ and $\fla$ from $0$.
And thus all the pre-images of $\al'(\la)$ under $P_\la^n$, $n\geq 0$ 
are separating points. 

The holomorphically moving fixed point $\al'(\la)$ for $f_\la$ 
has a unique pre-image $\al_1'(\la)\in\Klap$, 
which thus also moves holomorphically with $\la\in\WpqMz$. 
Moreover for $\la$ in any of the sub-limbs $\Lpqppqp$ 
the pre-image $\al_1'(\la)\in\Klap$ separates the critical value $-(\la/2)^2$ 
as well as $\al_1'(\la)$ itself from $\be':=0$. 
Define recursively $\al_n'(\la)$ as the unique $f_\la$ pre-image of $\al_{n-1}'(\la)$, 
which separates the later from $0$. Then $\al_n'\to 0$ as $n\to 0$ 
and the whole sequence $\{\al_n'(\la)\}$ moves holomorphically with 
$\la\in\WpqMz\Sm\rho^{-1}([0,1])$.  

\PROOF{of \thmref{dynamicsimplyparameter} :}\label{proofofdyn}
Let $\lastar\in\Mpq$ be any (Misiurewicz) parameter such that for some minimal $m\geq 1$ :  
$f_{\lastar}^m(-(\lastar/2)^2)=0$. 
Let $r>0$, $D_{\lambda^*}=\D(\lastar,r)$, and define 
the holomorphic function
$h(\la):=f_\la^m(-(\la/2)^2): D_{\lambda^*} \rightarrow~\C$, 
with $h(\lastar)=0$. 
Possibly restricting $r>0$ we can assume 
the sequence $\{\al_n'(\la)\}_n$ moves holomorphically  on  $D_{\lambda^*}$. 
By Slodkowskys Theorem there exists a holomorphic motion 
{\mapfromto H {D_{\lambda^*}\times\Chat} \C} of $\Chat$ 
extending the motion of the sequence $\{\al_n'(\la)\}_n$ 
union the $q$ external rays co-landing at $0$. 
We define a quasi-regular map {\mapfromto \zeta {D_{\lambda^*}} \C} by 
$\zeta(\la) = H_\la^{-1}(h(\la))$, where $H_\la$ 
is the quasi-conformal homeomorphism $H_\la(z) = H(\la,z)$. 
The map $\zeta$ is asymptotically conformal at $\lastar$, i.e. its
real dilatation $K_\zeta(\la)$ converge uniformly to $1$ as $\la\to\lastar$. 
More precisely its dilation at $\la$ is bounded by the dilatation of $H_\la$, 
and thus satisfies $\log K_\zeta(\la) \leq \log \d(\la)$, 
where $\d(\la)$ denotes the hyperbolic distance in $D_{\lambda^*}$ between $\lastar$ and $\la$. 

For each dynamical ray $R$ landing at $0$ for $\lastar$, 
each connected component of $\zeta^{-1}(R)$ is contained in a parameter ray 
landing on $\lastar\in\Mz$. 
It follows that the local degree of $\zeta$ at $\lastar$ is $1$, 
so that, possibly restricting $r>0$ further, the map $\zeta$ 
is a quasi-conformal homeomorphism from $D_{\lambda^*}$ 
into the dynamical plane of $f_\lastar$ such that 
for every $n$ with $\al_n'(\lastar)\in\zeta(D_{\lambda^*})$ and 
$\la_n := \zeta^{-1}(\al_n'(\lastar))$ we have 
$$
f_{\la_n}^{m+n}(-(\la_n/2)^2) = \al_n'(\la_n).
$$
Possibly reducing $r>0$ further we can also assume that 
$\zeta(D_{\lambda^*})$ is contained in a domain of univalence 
for the linearizing coordinate $\phi_\lastar$
for $f_\lastar$ at $0$. 
Let $N\in\N$ be such that $\al_n'(\lastar)\in\zeta(D_{\lambda^*})$ for all $n\geq N$. 
Then $\psipq := \phi_\lastar\circ\zeta$ is a quasi-conformal homeomorphism, 
which is asymptotically conformal at $\lastar$ and which maps the sequence 
$\{\la_n\}_{n\geq N}$ to the geometric sequence 
$\{\whal_n'(\lastar):=\phi_\lastar(\al_n'(\lastar))\}_{n\geq N}$ with 
$$
(\lastar)^q\cdot \whal_{n+1}'(\lastar) = \whal_n'(\lastar).
$$

For $\mustar = \xi(\lastar)$ we can construct similarly a 
quasi-conformal homeomorphism $\psiPQ$ 
defined on a small circular disk around $\mustar$ 
and asymptotically conformal at $\mustar$, mapping for large $n$ the sequence 
$\mu_n=\xi(\la_n)$ onto a geometric sequence 
$\{\whal_n'(\mustar)\}_{n\geq N}$ (possibly augmenting $N$) with 
$$
(\mustar)^Q\cdot \whal_{n+1}'(\mustar) = \whal_n'(\mustar).
$$
Let {\mapfromto {\whxi} {D_{\lambda^*}} \C} 
be any quasi-conformal extension of $\xi$ to $D_{\lambda^*}$, 
and consider the quasi-conformal map $\psi:=\psiPQ\circ\whxi\circ\psipq^{-1}$ 
between neighbourhoods of zero. 
It maps the $\lastar$ geometric sequence $\{\whal_n'(\lastar)\}_{n\geq N}$ 
onto the $\mustar$ geometric sequence $\{\whal_n'(\mustar)\}_{n\geq N}$. 
Moreover it maps the initial ray segment say $R^\lastar_{\theta_{p/q}}([0,\tau])$ 
onto a curve, which emanates from $0$ and which is homotopic relative to 
$\{0\}\cup\{\whal_n'(\mustar)\}_{n\geq N}$ 
(i.e. via a homotopy which fixes point wise this set and is injective above this set) 
to the ray segment $R^\mustar_{\theta_{P/Q}}([0,\tau])$.

We shall prove, \thmref{weakteichmullermap} below, that such a quasi-conformal homeomorphism satisfies 
$$
\limsup_{z \rightarrow 0} \Log K_{\psi}(z) \geq d_{\Hplus}(\Lastar, \Mstar).
$$
From this the Theorem follows. 
The following Theorems are generalizations of the {\Teichmuller} extremal theorem 
for complex tori, to a non-compact setting. 
\ENDPROOF
 \REFTHM{weakteichmullermap}
For $j = 1, 2$ let $\la_j\in\C\sm\Dbar$,   $\La_j$ be a logarithm of $\la_j$, and
{\mapfromto {\ga_j} {[0,1]} {\C}} be the arc $\ga_j(t) = -t^{\La_j}$.
Let {\mapfromto \psi D \C}, $D\supset\Dbar$ be a quasi-conformal homeomorphism  
with $\psi(\la_1^{-n}) = \la_2^{-n}$ for $n\geq 0$ 
and $\psi\circ\ga_1$ isotopic to $\ga_2$ relative to 
the set $\{0\}\cup\{\la_2^{-n}| n \geq 0\}$.
Then 
$$
\limsup_{z \rightarrow 0} \Log K_{\psi}(z) \geq d_{\Hplus}(\La_1, \La_2).
$$
\ENDTHM

\PROOF{of \thmref{weakteichmullermap} :}\label{proofofweak}
Notice that $\ga_j(t) = \exp(i\pi+\La_j\log t)$ for $0<t\leq 1$.
For each $j$ write $\whga_j^\pm(t) = \exp(\pm i\pi + t\La_j)$ for $t \in \R$. 
Then $\Log (\whga_j^-)$ is below while $\Log (\whga_j^+)$ is above  the bi-infinite sequences $\La_j\cdot\Z$, 
and together they separate $\La_j\cdot\Z$
from all the other pre-images $\log(\la_j^\Z)$ of $\la_j^\Z$.

Let {\mapfromto H {[0,1]\times[0,\tau]} {\C\sm\{\la_2^{-n}\}_{n\geq N}}} 
be a homotopy of $\psi(\ga_1([0,\tau]))$ to $\ga_2([0,\tau])$, 
with $H(s,0) = 0$, $H(s,t)\not=0$, for $t>0$, 
$H(0,t) = \psi(\ga_1)(t)$ and $H(1,t) = \ga_2(t)$. 
Choose $r'>0$ such that $D':=\D(0,r') \subset\subset D$. 
Possibly reducing $\tau$ we can suppose $H([0,1],[0,\tau])\subset D'$.
Possibly increasing $N$ we can suppose that $\la_1^{-n}\in D'$ for $n\geq N$ and 
for all $s\in[0,1]$ : $|H(s,\tau)| > |\la_2^N|$.
Possibly replacing $\psi(z)$ by $\la_2^N\cdot\psi(\la_1^{-N}z)$ 
we can suppose $\psi$ is quasi-conformal in a neighborhood $D'$ of $\Dbar$,  
$\psi(\la_1^{-n}) = \la_2^{-n}$ for all $n\geq 0$,  in particular $\psi(1)=1$ 
and $H([0,1],[0,\tau])\subset D'\sm\Dbar$. 
Let $\Psi$ be the unique lift of $\psi(\e^z)$ to $\e^z$ 
with $\Psi(0) = 0$. 
Then since $\Psi(0) = 0$ and $\psi(\ga_1([0,\tau]))$ 
is homotopic to $\ga_2([0,\tau])$ 
relative to $\{\la_2^n\}_{n\in\N}\cup\{0\}$ 
it follows that 
\REFEQN{pointconjugacy}
\forall\; n\geq 0, \forall\; m\in\Z :
\Psi(-n\La_1+m2\pi i) = -n\La_2 + m2\pi i.
\ENDEQN
Indeed let $\whHpm$ denote the lifts of $H$ to the covering $\e^z$ 
with $\whHpm(1,t) = \pm i\pi + t\La_j$. 
Then $\Re(\whHpm(s,\tau) > 0$ and since $\whHpm$ are homotopies, 
the curves $\Psi(\pm i\pi + ]-\infty,\tau]\cdot\La_1) = \whHpm(0,]-\infty,\tau])$ 
union $\whHpm([0,1],\tau)\cup(\pm i\pi +[-\infty,\tau[\cdot\La_2)$ 
are below for $-$ and above for $+$ the sequence $\{-n\La_2\}_{n\geq 0}$.\\ 
Thus \eqref{pointconjugacy} holds for $m=0$ and all $n\geq 0$.
And then \eqref{pointconjugacy} holds also for all $m$ since 
$\Psi(z+i2\pi) = \Psi(z) + i2\pi$.
The proof now follows from the following \thmref{WeakTeichmullermap}.
\ENDPROOF

\REFTHM{WeakTeichmullermap}
Suppose $R>0$, $\La_1, \La_2 \in \Hplus$, and {\mapfromto \Psi {\{x+iy| x < R \}} \C} 
is a quasi-conformal homeomorphism which satisfies 
\eqref{pointconjugacy}. Then 
$$
\limsup_{\Re(z)\to -\infty}\log K_\Psi(z) \geq \d_\Hplus(\La_1, \La_2).
$$
\ENDTHM

\PROOF
The following proof is based on a classical proof of the {\Teichmuller} extremal 
theorem. In this theorem the hypothesis on $\Psi$ above is replaced by 
a globally defined conjugacy {\mapfromto \Psi \C \C} between 
the full group actions in stead of only the semi-group actions on the semi orbit of $0$:
$$
\forall z\in\C, \forall n, m \in\Z : \Psi(z-n\La_1+m2\pi i) = \Psi(z)-n\La_2 + m2\pi i
$$
The proof of the {\Teichmuller} extremal theorem we lean on, 
amounts to showing the following :\\
\noindent\textbf{Claim :}
For any $\eps>0$ there exists integers $p, q, r, s$, $0<q$, $0\leq s$ 
such that the quadrilaterals $Q_j = Q(q\La_j-p2\pi i, s\La_j-r2\pi i)$ 
with corner $0$ and spanned by the vectors $q\La_j-p2\pi i$ and $s\La_j-r2\pi i$, 
say with $a$-sides parallel to the first vector have moduli satisfying 
$$
\log\frac{\Mod(Q_1)}{\Mod(Q_2)} \geq \d_\Hplus(\La_1, \La_2) - \eps.
$$
To turn the Claim into a proof in the case at hand with much weaker hypothesis 
we need to overcome that we have little control over the images of the boundary of $Q_1$. 
Let $\eps>0$ be given and choose integers $p, q, r, s$ as above with 
$q < 0$, $s \leq 0$ as above and consider the quadrilaterals 
$Q_j^k = -k\cdot Q_j\subset\Hminus$, $j = 1, 2$. 
Then $Q_j^k$ has the same modulus as $Q_j$, and $\Psi$ maps 
at least $k+1$ lattice points on each side of $Q_1^k$ 
to corresponding lattice points on $Q_2^k$. 
Since $\Psi$ is $K$-qc. for some fixed $K$ and conjugates the actions by the semi-groups 
on $0$ and since integer scaling maps the two orbits of $0$ into themselves. 
It follows that $\Psi(kz)/k$ converges locally uniformly to the affine 
mapping, which takes $\La_1$ to $\La_2$ and $i2\pi$ to itself. 
By continuity of the modulus it follows that we can choose $N$ such that 
for $k\geq N$ 
$$
\log\frac{\Mod(Q_2^k)}{\Mod(\Psi(Q_1^k)} \leq  \eps.
\quad\textrm{and thus}\quad 
\log\frac{\Mod(Q_1^k)}{\Mod(\Psi(Q_1^k)} \geq \d_\Hplus(\La_1, \La_2) - 2\eps.
$$
Since we can repeat this argument with uniform estimates on quadrilaterals to the 
left of any pair of corresponding lattice points the Theorem follows.

For completeness of exposition we include a proof of the claim above.

\noindent\textbf{Case I :} $\Im(\La_1) = \Im(\La_2) = 2\pi y$ 
and say $0<\Re(\La_2) < \Re(\La_1)$ 
so that $\d_\Hplus(\La_1,\La_2) = \log(\Re(\La_1)/\Re(\La_2))$. 
If $y=u/v$ with $(u,v)=1$ and $v>0$ take $p=u$, $q=-v$, $r=1$ and $s=0$.

If $y$ is irrational consider say the sequence of best rational approximants 
$u_n/v_n$ to $y$ with $v_n>0$. Then the quadrilaterals spanned by 
$v_n\La_j - u_n2\pi i$ and $2\pi i$ becomes more and more retangular 
with the same height $2\pi$ and ratio of the base $\Re(\La_1)/\Re(\La_2)$ 
so that given $\eps>0$ 
we may take $q=v_n$, $p=u_n$, $s=0$ and $r=1$ for a sufficiently large $n$.

\noindent\textbf{Case II :} $\Im(\La_1) \not= \Im(\La_2)$, 
say $\Im(\La_1) > \Im(\La_2)$ the other case being symmetric. 
Let $\Ga$ be the hyperbolic geodesic in $\Hplus$ through $\La_1$ and $\La_2$. 
Then $\Ga$ is a half circle orthogonal to $i\R$ 
with upper end point $i2\pi y$ for some $y\in\R$ and thus 
$2\pi y>\Im(\La_1)> \Im(\La_2)$.\\
\noindent\textbf{Case IIA :} 
If $y=u/v$ is an irreducible rational let $m/n$ be a rational with 
$nu-vm = 1$ and thus with $u/v > m/n$
Define 
$$
L_j := m\La_j - n2\pi i, \quad N_j := v\La_j - u2\pi i
\quad\textrm{and}\quad 
\whLa_j :=-i2\pi L_j/N_j, \quad j= 1, 2.
$$
and note that $-\pi/2 < \Arg(N_j) < 0$ for both values of $j$, 
as $i2\pi y$ is the upper end point of $\Ga$.
Then 
$$
\whLa_j = A(\La_j), \qquad\textrm{where}\qquad
A_m(z) := -i2\pi \frac{n z - m 2\pi i}{v z - u 2 \pi i} : \Hplus \to \Hplus
$$
(is conjugate to the {\Mobius} transformation 
$\frac{n z + m}{v z + u}\in\PSL(2,\Z)$ with determinant 
$nu-qv=1$ by the linear map $z\mapsto -i2\pi z$)
so that $\d_\Hplus(\La_1, \La_2) = \d_\Hplus(\whLa_1, \whLa_2)$ 
and $\Im(\whLa_1) = \Im(\whLa_2)$.
Thus we are essentially back to Case I. 
Let $\eps>0$ be given and let $a/b$ with $(a,b)=1$ satisfy 
$$
\log\frac{\Mod(Q(b\whLa_1-ai2\pi, -i2\pi)}{\Mod(Q(b\whLa_2-ai2\pi, -i2\pi)} 
\geq \d_\Hplus(\La_1, \La_2) - \eps.
$$
Then 
Possibly increasing $|a|$ and $b$ we can suppose 
that $-\pi/2 < \Arg((b\whLa_j-ai2\pi)\cdot N_j)< 0$, 
so that $(b\whLa_j-ai2\pi)N_j/(-i2\pi) = bL_j+aN_j\in\Hplus$ 
Thus we may take $(p,q,r,s) = (bm+au, bn+av, u, v)$.

\noindent\textbf{CASE IIB :} If $y$ is irrational. we may proceed similarly as above.
Let $u_m/v_m$, $v_m >0$ be the (odd numbered) convergents of $y$, 
which converges to $y$ from above, and let 
$u_{m+1}/v_{m+1}$ with $v_m>0$ be the following (even numbered) convergent 
Then $v_mu_{m+1} -u_mv_{m+1} = 1$, and the hyperbolic geodesic $\delta_m$ 
connecting $i2\pi u_m/u_m$ to 
$i2\pi u_{m+1}/v_{m+1}$ crosses $\Ga$ 
and converge to the upper end $i2\pi y$ of $\Ga$. 
In particular the crossing of $\Ga$ and $\delta_m$ takes place hyperbolically further 
and further away from the segment of $\Ga$ between $\La_1$ and $\La_2$.
Define 
$$
L_{j,m} := v_{m+1}\La_j - u_{m+1}2\pi i,\qquad N_{j,m} := v_m\La_j - u_m2\pi i 
\qquad\textrm{and}
$$
$\whLa_{j,m} :=-i2\pi L_{j,m}/N_{j,m}$.
Then $\Im(N_{j,m}) < 0$ 
and $\Arg(\whLa_{1,m}-\whLa_{2,m})$ converge to $0$. 
Indeed, the hyperbolic isomorphism
$$
A_m(z) := -i2\pi \frac{v_{m+1} z - u_{m+1} 2\pi i}{u_m z - u_m 2 \pi i}
$$
maps $i2\pi u_m/v_m$ to $\infty$, $i2\pi u_{m+1}/v_{m+1}$ to $0$ 
and $\delta_m$ to $\Rplus$. 
So the hyperbolic isomorphism $A_m$ also maps the geodesic $\Ga$ 
to a semi-circle geodesic $\whGa$ crossing $\Rplus$, and $\La_1, \La_2$ 
to $\whGa$ hyperbolically further and further away from $\Rplus$. 
Moreover as $u_m/v_m-y > y-u_{m+1}/v_{m+1}$ the angle between 
$\delta$ and $\Ga$, when these are parametrized with increasing imaginary parts, is 
smaller than $\pi/2$ for sufficiently large $m$.
This proves $\Arg(\whLa_{1,m}-\whLa_{2,m})$ converge to $0$. 
Let $a_m/b_m$ be the irreducible rational with smallest denominator between 
$\Im(\whLa_1)/2\pi$ and $\Im(\whLa_2)/2\pi$.
Given any $\eps>0$ we may choose a sufficiently large 
$m$ such that 
$$
\log\frac{\Mod(Q(b\whLa_1-ai2\pi, -i2\pi)}{\Mod(Q(b\whLa_2-ai2\pi, -i2\pi)} 
\geq \d_\Hplus(\La_1, \La_2) - \eps.
$$
From here proceed as in the CASE IIB, $y$ rational above to produce the 
desired $(p,q,r,s)$. 
\ENDPROOF

\subsection{Multiplier relation at simple satelite parabolic bifurcations}\label{relatingmultipliers}
For {\mapfromto F U \C}, where $U\subset\C$, a holomorphic map 
with a fixed point $z_0\in U$, we define:
\begin{align}\label{resitetc}
\mult(F,z_0) :&= \frac{1}{2\pi i} \cint_{C(z_0,r)} \frac{F'(w) - 1}{F(w)-w}\dw,\\
&=\textrm{the multiplicity of $z_0$ as a fixed point for $F$.}\nonumber\\
\iota(F,z_0) :&= \frac{1}{2\pi i} \cint_{C(z_0,r)} \frac{1}{w-F(w)}\dw,\\
&= \textrm{the holomorphic index of $F$ at $z_0$.}\nonumber\\
\resit(F,z_0) :&= \frac{1}{2}\mult(F,z_0) - \iota(F,z_0),\\
&=\textrm{the \emph{\residuiteratif} of $F$ at $z_0$.}\nonumber
\end{align}
Here and elsewhere the radius of the circle for the circulation is chosen 
sufficiently small for $z_0$ to be the only fixed point of $F$ inside the circle, 
except if explicitly stated otherwise. 
The {\residuiteratif} is well behaved under iteration: indeed
for all $n\geq 1$ we have $n\cdot\resit(F^n,z_0) = \resit(F,z_0)$
(see e.g. \cite[p. 250]{BuffandEpstein}).
For a domain $U$ on which the principal logarithm $\Log F'$ is well defined and holomorphic 
we define the Buff-form
$$
\om_F := \frac{(F'(z) - 1)\dz}{(F(z)-z)\Log F'(z)}.
$$
A routine calculation shows that, for a fixed point $z_0$ as above:
\REFEQN{logresidue}
\frac{1}{2\pi i} \cint_{C(z_0,r)} \frac{(F'(w) - 1)\dw}{(F(w)-w)\Log F'(z)} = 
\begin{cases}
\frac{1}{\Log F'(z_0)} & \textrm{if } F'(z_0)\not=1,\\
\resit(F,z_0) & \textrm{if } F'(z_0) = 1
\end{cases}
\ENDEQN
(see also the habilitation thesis of Cheritat \cite{Cheritat}). 

Suppose that, for some irreducible rational $p/q$ and $\opq=\e^{i2\pi p/q}$, the family:
$$
F_\la(z) = \la z + \OO(z^2), \quad\textrm{with}\quad 
\la \in\D(\opq,R), \quad R>0
$$
is a holomorphic family of holomorphic maps defined on a neighourhood of $0$ 
and with first $q$ iterates defined and holomorphic on some disk $\D(r'),\, r'>0$. 
Suppose furthermore that $0$ is a non degenerate parabolic fixed point 
for $F=F_\opq$, \textit{i.e.} $\mult(F^q,0) = q + 1$. 
Choose $r,~0< r < r'$, such that $0$ is the unique fixed point of $F$ 
in $\overline{\D(r)}$, and for every $z \in\overline{\D(r)}$,  $(F^q)'(z)\in\D(1,1)$. 
Then, possibly reducing $R>0$ and $r\in]0,r'[$, 
we can assume that, for all $\la\in\D(\opq,R)$ and for all $z \in \D(r)$,
$(F_\la^q)'(z)\in\D(1,1)$, and the circle $C(0,r)= \partial \D(r)$
is disjoint from the set of fixed points of $F_\la^q$.
For such $\la$ the disk $\D(r)$ contains precisely 
$q+1$ fixed points of $F_\la^q$, and no other periodic point of $F_\la$ 
with a period dividing $q$. 
Moreover, for all $\la$ the fixed point $0$ is non degenerate, 
and for $\la\not=\opq$ the other $q$ fixed points of $F_\la^q$ form a 
non degenerate $q-cycle$ for $F_\la$. 
We denote by $\La$ the principal logarithm of the multiplier $\la^q$ of the fixed point $0$ for $F_\la^q$. 
We denote by $\rho=\rho(\la)$ the multiplier of the $q$-cycle for $F_\la$, and by 
\emph{$\mathcal{P} =\Log\rho$} the principal logarithm of $\rho$. 
Then, by the hypothesis above, both $\La$ and $\PPP$ belong to the disk 
$\Log(\D(1,1))$, and we may take $\La$ as the parameter. 
We shall abuse the notation and write $F_\La$ for $F_\la$ and 
$\PPP(\La)$ for $\PPP(\la)$. Define 
$$
H(\La) = \frac{1}{2\pi i} \cint_{C(0,r)} \frac{((F_\La^q)'(w) - 1)\dw}
{(F_\La^q(w)-w)\Log (F_\La^q)'(z)},
$$ 
where $C(0,r)$ is the circle discussed in the paragraph immediately above. 
Then by \eqref{logresidue} the holomorphic function $H$ satisfies 
\REFEQN{Resitfunction}
H(\La) = 
\begin{cases}
\resit(F_0^q,0), &\textrm{if}\quad \La=0,\\
\frac{1}{\La} + \frac{q}{\PPP(\La)}, &\textrm{if}\quad \La\not=0.
\end{cases} 
\ENDEQN
When $\La\not=0$ we obtain, by solving the above equation for $\PPP(\La)$, 
the formula
$$
\PPP(\La) = -\frac{q\La}{1 - \La H(\La)}
$$
which evidently is also valid for $\La=0$. It shows that 
$\PPP'(0) = -q \not=0$. Thus, by the inverse function theorem and possibly 
at the expense of reducing $R>0$ further, we may take $\PPP$ as the parameter 
and reparametrize again. We shall thus write $\La(\PPP)$ rather than 
$\PPP(\La)$, and $\Res_{p/q} := \resit(F_0^q,0)$.
\REFPROP{LambdaofP}
The holomorphic function $\La(\PPP)$ has the following second order Taylor expansion 
$$
\La(\PPP) = -\PPP/q - \Res_{p/q}(\PPP/q)^2 +  \OO((\PPP/q)^3).
$$
\ENDPROP
\PROOF
Solving this time the second line of \eqref{Resitfunction} for $\La(\PPP)$ we find 
$$
\La(\PPP) = -\frac{\PPP}{q - \PPP H(\La(\PPP))} = 
-\frac{\PPP}{q}\left(1 + \frac{\PPP H(\La(\PPP))}{q} + 
\OO\left(\left(\frac{\PPP H(\La(\PPP))}{q}\right)^2\right)\right)
$$
As $H$ is holomorphic with $H(0) = \Res_{p/q}$ we may write 
$$
H(\La(\PPP)) = \Res_{p/q} + \Res_1\La(\PPP) + \OO(\La^2(\PPP))
= \Res_{p/q} - \frac{\PPP\Res_1}{q} + \OO\left(\left(\frac{\PPP}{q}\right)^2\right),
$$ 
where $\Res_1= H'(0)$. Inserting this into the 
formula for $\La(\PPP)$ above we obtain the Proposition. 
\ENDPROOF

Buff and Epstein \cite{BuffandEpstein} have proved inequalities for the {\residuiteratif}. 
Recall that $P_{\lambda}(z)=\lambda z+z^2$.
As a particular instance of the immediate Corollary of their Theorem B we have:
\REFTHM{BuffEpstein}
For any irreducible rational $p/q$ and $\la_0= \e^{i2\pi p/q}$ 
the {\residuiteratif} of $\P_{\la_0}^q$ at $0$ satisfies
$$
\Re(\resit(P_{\la_0}^q,0)) \geq 1.
$$
\ENDTHM
\PROOF
For $p/q = 0/1$ an explicit computation for the polynomial $P_1(z)=z+z^2$ 
yields $\resit(P_1,0) =1$. For $q>1$ we have from Buff and Epsteins formula 
\cite[Corollary of Theorem B p.~252]{BuffandEpstein}:
$$
\Re(\resit(P_{\la_0}^q,0)) \geq \frac{1}{2\log 2} + \frac{q}{4} > 1.
$$
\ENDPROOF

\PROOF{\textit{of \propref{goodassymtoticestimates}: }\label{proofofassympdevelop}}
For the asymptotic formulas, take $F_\la = \Pla$ and $F_\mu = \Pmu$ 
in the above discussion and \propref{LambdaofP}. 
For the lower bound on the real part of $\Res_{p/q}= \resit(F_0^q,0)$, apply 
\thmref{BuffEpstein}.
\ENDPROOF

\PROOF{\textit{of \corref{multiplierhyperbolicdistancediverge}: }}
Recall that, for $\rho = \e^{it}\in\Dbar$, we have $M(\rho) = \xi(\La(\rho))$. 
We compute the hyperbolic distance between $M(\rho)$ and $\La(\rho)$ from 
the asymptotic formula given by \propref{goodassymtoticestimates}. 
To this end note that, since the real parts of both $R_1:=\Res_{p/q}$ and 
$R_2:=\Res_{P/Q}$ are positive, the second term of the asymptotic formulas:
$$
- \Res_{p/q}\left(\frac{\PPP}{q}\right)^2=-R_1\left(\frac{it}{q}\right)^2 = R_1\left(\frac{t}{q}\right)^2
\quad\textrm{and}\quad $$
$$
- \Res_{P/Q}\left(\frac{\PPP}{Q}\right)^2=-R_2\left(\frac{it}{Q}\right)^2 = R_2\left(\frac{t}{Q}\right)^2
$$
have positive real parts, and therefore $\La(\rho)$ and $M(\rho)$ 
both belong to $\Hplus$. 
The following computation uses the fact that homothetisies and vertical translations are hyperbolic 
isometries in $\Hplus$.
\ALIGN
\d_\Hplus(\La(\rho),M(\rho)) &= 
\d_\Hplus(-\frac{it}{q}+\frac{t^2R_1}{q^2} + \OO((t/q)^3), 
-\frac{it}{Q}+\frac{t^2R_2}{Q^2} + \OO((t/Q)^3))\\
(\frac{Qq}{t^2}\cdot~)\qquad&=\d_\Hplus(-\frac{iQ}{t}+R_1\frac{Q}{q}+ \OO\left(\frac{tQ}{q^2}\right),
-\frac{iq}{t}+R_2\frac{q}{Q}+ \OO\left(\frac{tq}{Q^2}\right))\\
&\geq \d_\Hplus(-\frac{iQ}{t}+R_1\frac{Q}{q},
-\frac{iq}{t}+R_2\frac{q}{Q}) - \OO(t)\\
&\geq \d_\Hplus(-\frac{iQ}{t}+1,
-\frac{iq}{t}+1) - \OO(1).
\end{align*}
Finally, hyperbolic geodesics in $\Hplus$ 
are arcs of circles with center on the imaginary axis, 
and with diameter their intersection with the imaginary axis. 
Since the chord of a circle is smaller than the diameter, we obtain: 
\REFEQN{logboundonhypdist}
\d_\Hplus(\La(\rho),M(\rho)) \geq 2\log\frac{Q-q}{t} - \OO(1)\longrightarrow \infty,
\qquad\textrm{as}\qquad t\rightarrow 0.
\ENDEQN
\ENDPROOF

\subsection{Yoccoz estimates of quadratic limbs.}
\setcounter{theorem}{6}
\PROP
For any irreducible rational $p/q$ the $\frac{n^2-1}{n^3}$-sublimbs of $\Lpq$ 
satisfy
$$
\diam_\Hplus\left(\whLpqn\right),~
\longrightarrow~ 0\qquad\textrm{as}\quad n\rightarrow \infty.
$$
\ENDPROP

\PROOF{}
Note at first that, for every integer $n\geq 2$, 
the rational number $(n^2-1)/n^3$ is irreducible and asymptotic to $1/n$. 
Recall that $\whHpq$ with root $0$ is the principal hyperbolic component of $\whMpq$. 
Let $ C_{p/q}>0$ be the Yoccoz-constant for $\whHpq$  
given by  \eqref{yoccozparameterineq} so that 
\REFEQN{Euclideansmalllimbs}
\diam(\whLpqn) \leq \frac{C_{p/q}}{n^3}.
\ENDEQN
Moreover, by \propref{goodassymtoticestimates}, the root 
$\La(\e^{i2\pi \frac{n^2-1}{n^3}})$ 
of the limb $L_{\frac{n^2-1}{n^3}}^{p/q}$ is placed at 
$$
\La(\e^{i2\pi \frac{n^2-1}{n^3}}) = -\frac{i2\pi \frac{n^2-1}{n^3}}{q} + 
\Res_{p/q}\left(\frac{2\pi \frac{n^2-1}{n^3}}{q}\right)^2 
+ \OO\left(\left(\frac{2\pi \frac{n^2-1}{n^3}}{q}\right)^3\right),
$$
so that
\REFEQN{rootdistancetobdry}
\Re(\La(\e^{i2\pi \frac{n^2-1}{n^3}})) \geq 
\Re(\Res_{p/q})\frac{4\pi^2}{q^2n^2} - \OO(1/n^3)
\ENDEQN
Combining the estimates \eqref{Euclideansmalllimbs} 
and \eqref{rootdistancetobdry} we find that 
\REFEQN{smalllimbsareoneovern}
\diam_\Hplus(\Log(L_{\frac{n^2-1}{n^3}}^{p/q})) \leq \OO(1/n)
\ENDEQN
\ENDPROOF

\PROOF{of \corref{Dynamicaldivergenceinsmalllimbs} :}
Let $\La_n$ and $M_n$ be the principal logarithms of the $\beta'$ fixed points 
for $f_{\lambda_n}$ and $g_{\mu_n}$ respectively, where $\mu_n=\xi(\lambda_n)$ and $n\geq 2$. Define 
$\La_n^0 = \La(\rho_n)$ and $M_n^0 = M(\rho_n)$, where 
$\rho_n = \e^{i2\pi{\frac{n^2-1}{n^3}}}$ (so that 
$M_n^0 = \xi(\La_n^0)$).
Then by the triangle inequality:
$$
\d_\Hplus(\La_n,M_n) 
\geq \d_\Hplus(\La_n^0,M_n^0) - \d_\Hplus(\La_n^0,\La_n) - \d_\Hplus(M_n,M_n^0).
$$
The first term on the righthandside diverge to $\infty$ by \corref{multiplierhyperbolicdistancediverge}, 
and the two remaining terms converge to $0$ by \propref{hyperbolicallysmalllimbs}. 
In fact asymptotically we have the more explicit lower bound
\eqref{logboundonhypdist} for the first term and the upper bound 
\eqref{smalllimbsareoneovern} for the two remaining terms. 

\ENDPROOF

\end{document}